\documentclass{article}
\usepackage{amsmath,amsthm,amssymb,authblk}

\newtheorem{theorem}{Theorem}[section]
\newtheorem{corollary}[theorem]{Corollary}
\newtheorem{proposition}[theorem]{Proposition}
\newtheorem{lemma}[theorem]{Lemma}

\numberwithin{equation}{section}
\def\pf{\noindent{\bf Proof.}~}

\begin{document}
\title{Permutation-like Matrix Groups with a Maximal Cycle of Power of
Odd Prime Length }
\author{Guodong Deng,\quad Yun Fan\\
\small School of Mathematics and Statistics\\
\small Central China Normal University, Wuhan 430079, China}
\date{}
\maketitle

\insert\footins{\footnotesize{\it Email address}:
yfan@mail.ccnu.edu.cn (Yun Fan).}

\begin{abstract}
If every element of a matrix group is similar to a permutation matrix,
then it is called a permutation-like matrix group. References \cite{C07}
and \cite{DF} showed that, if a permutation-like matrix group
contains a maximal cycle of length equal to a prime or a square of a prime
and the maximal cycle generates a normal subgroup,
then it is similar to a permutation matrix group.
In this paper, we prove that if a permutation-like matrix group
contains a maximal cycle of length equal to any power of any odd prime
and the maximal cycle generates a normal subgroup,
then it is similar to a permutation matrix group.


\end{abstract}

\section{Introduction}

Let ${\Bbb C}$ be the complex filed. Let
${\rm GL}_d({\Bbb C})$ be the complex general linear group
of dimension $d$, i.e. the multiplicative group consisting of
invertible complex $d\times d$ matrices.
Any subgroup ${\cal G}\le{\rm GL}_d({\Bbb C})$ is called a matrix group of
dimension $d$. For a matrix group ${\cal G}$, if there exists a
$P\in{\rm GL}_d({\Bbb C})$ such that $P^{-1}AP$ for all $A\in{\cal G}$
are permutation matrices, then ${\cal G}$ is said to be
{\em similar (or conjugate) to a permutation matrix group},
or said to be a {\em permutation matrix group} for short.
If for every element $A$ of ${\cal G}$ there exists an
$S\in{\rm GL}_d({\Bbb C})$ such that $S^{-1}AS$ is a permutation matrix, 
then ${\cal G}$ is called a {\em permutation-like matrix group}.

Cigler \cite{C05,C07} showed that a permutation-like matrix group
is not a permutation matrix group in general.
A $d\times d$ matrix is called a {\em maximal cycle} if
it is similar to a permutation matrix corresponding
to the cycle permutation of length~$d$.
Cigler conjectured that:

\hangindent13ex
\noindent{\it Conjecture:}~
{If a permutation-like matrix group contains a maximal cycle, then
it is a permutation matrix group.}

Cigler \cite{C05,C07} proved that, if a permutation-like matrix group ${\cal G}$
of dimension prime $p$ contains a maximal cycle which generates
a normal cyclic subgroup, then ${\cal G}$ is a permutation matrix group.
In \cite{DF}, we further proved that this result still holds 
if we replace the dimension prime $p$ by dimension $p^2$.

In this paper we prove:

\begin{theorem}\label{main}
Let ${\cal G}$ be a permutation-like matrix group of dimension $p^n$ 
where $p$ is any odd prime and $n$ is any positive integer.
If ${\cal G}$ contains a maximal cycle which generates
a normal cyclic subgroup, then ${\cal G}$ is a permutation matrix group.
\end{theorem}

Necessary preparations for the proof of the theorem are made in Section \ref{pre}.
For fundamentals of the group theory, please refer to \cite{AB,R}.
In Section \ref{case p group} we treat a special case of the theorem
where ${\cal G}$ is a $p$-group, i.e. its order is a power of $p$.
The proof of Theorem \ref{main} is presented in Section \ref{proof}.

\section{Preparation}\label{pre}

We begin with few preliminaries on relations between
a cyclic group $\langle C\rangle$ of order $d$ generated by $C$
and the residue ring ${\Bbb Z}_{d}$ of the integer ring ${\Bbb Z}$
modulo the positive integer $d$.
We have a bijection which maps $k\in{\Bbb Z}_{d}$ to $C^k\in\langle C\rangle$.
By ${\Bbb Z}_d^*$ we denote the multiplicative group consisting
of the reduced residue classes in ${\Bbb Z}_{d}$.  Obviously,
$\langle C^k\rangle=\langle C\rangle$ if and only if
$k\in{\Bbb Z}_d^*$; at that case we say that $C^k$
is a {\em generator} of the cyclic group $\langle C\rangle$.
Any automorphism $\alpha$ of $\langle C\rangle$ is corresponding to exactly
one $r\in{\Bbb Z}_d^*$ such that
$\alpha(C^k)=C^{rk}$ for any $C^k\in\langle C\rangle$,
By $\mu_r$ we denote the permutation of ${\Bbb Z}_d$
mapping $k\in{\Bbb Z}_{d}$ to $\mu_r(k)=rk$.
The automorphism group of $\langle C\rangle$
is isomorphic to ${\Bbb Z}_d^*$ by mapping $\alpha$ to $\mu_r$ as above.
If ${\cal G}$ is a finite group which contains $\langle C\rangle$
as a normal subgroup,
then ${\cal G}$ is homomorphic to a subgroup
of ${\Bbb Z}_d^*$ with kernel consisting of
the elements which centralize $\langle C\rangle$.

For an element $A$ of a group ${\cal G}$, we denote the order of $A$
by ${\rm ord}_{\cal G}(A)$, or ${\rm ord}(A)$ for short 
if the group is known from context.
Let $p$ be a prime. A finite group ${\cal G}$ is called a {\em $p$-group}
if its order $|{\cal G}|$ is a power of $p$.
The next lemma about a cyclic $p$-group is obvious.
As mentioned above, it is also a lemma about ${\Bbb Z}_{p^n}$.
Note that $|{\Bbb Z}_{p^n}^*|=p^{n-1}(p-1)$.

\begin{lemma} \label{p-cyclic}
Let ${\cal G}=\langle C\rangle$ be a cyclic group and
$|{\cal G}|=p^n$.

(i)~ For $0\le a<n$ denote ${\cal G}^{p^a}=\{C^{p^a t}\mid 0\le t<p^{n-a}\}$.
Then ${\cal G}^{p^a}=\langle C^{p^a}\rangle$
is a cyclic subgroup of ${\cal G}$ of order $p^{n-a}$ generated by $C^{p^a}$.

(ii)~ The map ${\cal G}\to{\cal G}^{p^a}$, $C^k\mapsto C^{p^a k}$,
is an epimorphism with kernel ${\cal G}^{p^{n-a}}$.
In particular,  for any generator $G$ of the cyclic group ${\cal G}^{p^a}$
there are exactly $p^{a}$ generators of ${\cal G}$ which are mapped to $G$.

(iii)~ If $p$ is an odd prime and $\alpha$ is an automorphism of ${\cal G}$,
then ${\rm ord}(\alpha)=sp^a$ with $s\mid (p-1)$ and $0\le a<n$,
and there are integers $u,v$ coprime to $p$ such that
${\rm ord}_{{\Bbb Z}_{p^n}^*}(u)=s$ 
(hence ${\rm ord}_{{\Bbb Z}_{p}^*}(u)=s$)
and $\alpha(C)=C^{u+vp^{n-a}}$.
\qed
\end{lemma}

For a subgroup ${\cal H}$ of a group ${\cal G}$,
if any element which centralizes ${\cal H}$ is contained in
${\cal H}$, then ${\cal H}$ is said to be {\em self-centralized}.

\begin{lemma} \label{split}
Assume that ${\cal G}$ is a finite group
containing a self-centralized normal cyclic $p$-subgroup $\langle C\rangle$
generated by an element $C$, where $p$ is an odd prime.
Then ${\cal G}=\langle A,C\rangle$ generated by two elements, and
one of the following two holds.

(i)~ If ${\cal G}$ is a $p$-group, then there is an
$A'\in{\cal G}$ such that ${\cal G}=\langle A',C\rangle$ and
$\langle A'\rangle\cap\langle C\rangle=1$.

(ii)~ If ${\cal G}$ is not a $p$-group, then
$\langle A\rangle\cap\langle C\rangle=1$.
\end{lemma}

\pf We sketch a proof, some of the arguments will be quoted later. 

Assume that $|\langle C\rangle|=p^n$.
Since $\langle C\rangle$ is self-centralized,
the quotient group ${\cal G}/\langle C\rangle
\cong\langle A\rangle/\langle A\rangle\cap\langle C\rangle$
is isomorphic to a subgroup of the multiplicative ${\Bbb Z}_{p^n}^*$.
The group ${\Bbb Z}_{p^n}^*$ is a cyclic group of order $p^{n-1}(p-1)$.
Thus, $|{\cal G}/\langle C\rangle|=sp^a$ with $0\le a<n$ and $s\,|\,(p-1)$,
and there is an $r\in{\Bbb Z}_{p^n}^*$ such that
\begin{equation}\label{r}
 A^{-1}CA=C^r~~\mbox{and}~~
 {\rm ord}_{{\Bbb Z}_{p^n}^*}(r)=|{\cal G}/\langle C\rangle|=sp^a.
\end{equation}
For any integer $k$ and positive integer $j$, it is a direct computation that
\begin{equation}\label{(AC)^j}
(AC^k)^j=(AC^k)\cdots(AC^k)(AC^k)=A^{j}C^{k(r^{j-1}+\cdots+r+1)}.
\end{equation}

(i).~ Assume that $|{\cal G}/\langle C\rangle|=p^a$,  $0<a<n$.
By Lemma \ref{p-cyclic}(iii), we assume that $r=1+vp^{n-a}$ with $p\nmid v$.
Denote $\langle C\rangle^A=\{C^t\mid A^{-1}C^tA=C^t\}$.
Note that $A^{-1}C^tA=C^{tr}$. So
$A^{-1}C^tA=C^t$ if and only if $t(r-1)\equiv 0\pmod{p^n}$.
Thus $\langle C\rangle^A =\langle C^{p^a}\rangle$.
Since ${\cal G}/\langle C\rangle
\cong\langle A\rangle/\langle A\rangle\cap\langle C\rangle$,
$A^{p^a}\in\langle C\rangle$ hence
$A^{p^a}\in\langle C\rangle^A=\langle C^{p^a}\rangle$.
So, we can find an integer $k$ such that $A^{p^a}C^{kp^a}=1$.
Note that
\begin{equation}\label{r-1}
r^{p^a-1}+\cdots+r+1=\frac{r^{p^a}-1}{r-1}
=\frac{(1+vp^{n-a})^{p^a}-1}{vp^{n-a}} \equiv p^a\pmod{p^n}.
\end{equation}
Let $A'=AC^k$. By Eqns \eqref{(AC)^j} and \eqref{r-1},  we obtain
$$
A'^{p^a}=A^{p^a}C^{k(r^{p^a-1}+\cdots+r+1)}=A^{p^a}C^{kp^a}=1.
$$
Thus, ${\cal G}=\langle A',C\rangle$ and
$\langle A'\rangle\cap\langle C\rangle=1$.

(ii).~ Assume that $|{\cal G}/\langle C\rangle|=sp^a$, $0\le a<n$ and $s\mid(p-1)$.
Since ${\cal G}$ is not a $p$-group, we have $s>1$.
By Lemma \ref{p-cyclic} (iii), we can assume that
${\rm ord}_{{\Bbb Z}_{p}^*}(r)=s$. Then $p\nmid(r-1)$ (otherwise $s=1$).
As we have seen, $A^{-1}C^tA=C^t$ if and only if $t(r-1)\equiv 0\pmod{p^n}$.
At the present case, $t(r-1)\equiv 0\pmod{p^n}$ 
if and only if $t\equiv 0\pmod{p^n}$. Thus $\langle C\rangle^A=1$.
As argued above, $A^{sp^a}\in\langle C\rangle^A$.
So $A^{sp^a}=1$. That is, $\langle A\rangle\cap\langle C\rangle=1$.
\qed

We turn to preliminaries on matrices.

A diagonal blocked matrix $\begin{pmatrix}B_1\\ & \ddots\\ &&B_m\end{pmatrix}$
is denoted by $B_1\oplus\cdots\oplus B_m$ for short.
The identity matrix of dimension $d$ is denoted by $I_{d\times d}$, or
$I$ for short if the dimension is known from context.

We denote the characteristic polynomial of a
complex matrix $A$ by ${\rm char}_A(x)$.

\begin{lemma}[{\cite[Lemma 2.1]{DF}}]\label{permutable}
The following two are equivalent to each other:

(i)~  $A$ is similar to a permutation matrix;

(ii)~ $A$ is diagonalizable and ${\rm char}_A(x)=\prod_{i}(x^{\ell_i}-1)$.

\noindent
If it is the case, then each factor $x^{\ell_i}-1$ of ${\rm char}_A(x)$
corresponds to exactly one $\ell_i$-cycle of the cycle decomposition
of the permutation of the permutation matrix. \qed
\end{lemma}

By $\Phi_n(x)$ we denote the cyclotomic polynomial of degree $n$, i.e.
$\Phi_n(x)=\prod_{\omega}(x-\omega)$ with $\omega$ running over
the primitive $n$-th roots of unity. 
Since $x^n-1=\prod_{k|n}\Phi_k(x)$, 
the following is an immediate consequence of the above lemma.

\medskip\noindent{\bf Corollary.}~
{\it Let $A$ be a matrix similar to a permutation matrix, and
$m,n$ be positive integers. If ${\Phi_n(x)}^m\,\big|\,{\rm char}_A(x)$,
then $\Phi_k(x)^m\,\big|\,{\rm char}_A(x)$ for any $k\,|\,n$.
\qed
}

\medskip
The next lemma is a combination of
\cite[Eqns (2.1), (2.2), (2.3) and (2.4)]{DF}.

\begin{lemma}\label{spectrum}
Let $C\in{\rm GL}_d({\Bbb C})$ be a maximal cycle of dimension $d$,
and $\lambda$ be a primitive $d$-th root of unity.
Then the following hold.

(i)~ $\{\lambda^j\mid j\in{\Bbb Z}_{d}\}
 =\{\lambda^0=1,\lambda,\cdots,\lambda^{d-1}\}$
is the spectrum (i.e. the set of eigenvalues) of $C$.

(ii)~ The eigen-subspace of every eigenvalue $\lambda^j$ of $C$,
denoted by ${\rm E}(\lambda^j)$, has dimension $1$, and
${\Bbb C}^{d}=\bigoplus_{j=0}^{d-1}{\rm E}(\lambda^j)$.

(iii)~ If $e_j\in {\rm E}(\lambda^j)$ for $j=0,1,\cdots, d-1$ are all nonzero, then
$(e_0,e_1,\cdots,e_{d-1})$ is a basis of ${\Bbb C}^{d}$ and
$Ce_j=\lambda^{j}e_j$, $j=0,1,\cdots,d-1.$

(iv)~ Let $f=\alpha_0e_0+\alpha_1e_1+\cdots+\alpha_{d-1}e_{d-1}$
where $e_0,\cdots,e_{d-1}$ are as above and all $\alpha_j\in{\Bbb C}$.
Then $(f, ~ Cf,~\cdots,~C^{d-1}f)$ is a basis of ${\Bbb C}^d$ if and only if
$\alpha_j\ne 0$ for all $j=0,1,\cdots,{d-1}$. 
If it is the case, with respect to that basis, 
$C$ is a cycle permutation matrix of dimension $d$. \qed
\end{lemma}

\begin{lemma}[{\cite[Lemma 2.2]{DF}}]\label{basis}
Let $C$ and $\lambda$ be as above in Lemma \ref{spectrum}.
Let $A\in{\rm GL}_d({\Bbb C})$ such that $A^{-1}CA=C^r$
for an $r\in{\Bbb Z}_d^*$.  Further assume that
${\Bbb Z}_d$ is partitioned into $\mu_r$-orbits as follows:
$$\Gamma_0=\{0\},~ \Gamma_1=\{j_1,rj_1,\cdots,r^{d_1-1}j_1\},~
 \cdots,~ \Gamma_m=\{j_m,rj_m,\cdots,r^{d_m-1}j_m\}$$
i.e. $r^{d_k}j_k\equiv j_k\pmod d$ but $r^{d_k-1}j_k\not\equiv j_k\pmod d$.
For $k=0,\cdots, m$, take nonzero $e_{j_k}\in {\rm E}(\lambda^{j_k})$, set
${\cal E}_k=\{e_{j_k},Ae_{j_k},\cdots,A^{d_k-1}e_{j_k}\}$
and ${\cal E}=\bigcup_{k=0}^m{\cal E}_k$. Then the following hold.

(i)~ $A\cdot{\rm E}(\lambda^{j})={\rm E}(\lambda^{rj})$, $j=0,\cdots,d-1$,
where $A\cdot{\rm E}(\lambda^{j})=\{Av_j| v_j\in{\rm E}(\lambda^{j})\}$.

(ii)~ ${\cal E}_k$ is a basis of
$V_k:=\bigoplus_{h=0}^{d_j-1}{\rm E}(\lambda^{r^hj_k})$,
and $A$ restricted to $V_k$ has the matrix
\begin{equation}\label{monomial cycle}
   A|_{{\cal E}_k}= \begin{pmatrix}0&\cdots&0&\omega_k\\ 1&\ddots&\ddots&0\\
    \vdots&\ddots&\ddots&\vdots\\ 0&\cdots&1&0\end{pmatrix}_{d_k\times d_k},
\end{equation}
where $\omega_k$ is an $({\rm ord}(A)/d_k)$-th root of unity,

(iii)~ ${\Bbb C}^d=V_1\oplus\cdots\oplus V_m$,
the union ${\cal E}={\cal E}_1\cup\cdots\cup{\cal E}_m$ is a basis
of ${\Bbb C}^d$ and, with respect to the basis ${\cal E}$, the matrix of $A$ is
\begin{equation*}
 A|_{\cal E}=A|_{{\cal E}_1}\oplus\cdots\oplus A|_{{\cal E}_m}.
 \eqno\qed
 \end{equation*}
\end{lemma}

\begin{proposition}[{\cite[Proposition 2.3]{DF}}]\label{SC}
Let the notations be as in Lemma \ref{basis}.
If the matrix $A|_{\cal E}$ is a permutation matrix (i.e. all $\omega_k=1$),
then the matrix group $\langle A,C\rangle$ generated by $A$ and $C$
is a permutation matrix group. \qed
\end{proposition}

Note that ``the matrix $A|_{\cal E}$ is a permutation matrix''
means that mapping $e\in{\cal E}$ to $Ae$ is a permutation of the set ${\cal E}$.

\begin{lemma}[{\cite[Proposition 4.2]{C07}}]\label{self-centralized}
If ${\cal G}$ is a permutation-like matrix group and
$C\in{\cal G}$ is a maximal cycle, then $\langle C\rangle$ 
is self-centralized in ${\cal G}$. 
\end{lemma}

We sketch the proof for convenience. For $A\in{\cal G}$ which centralizes $C$, 
we can assume that $C$ is a permutation matrix of dimension $d$ corresponding to 
the $d$-cycle permutation and $A=\sum_{i=0}^{d-1}\alpha_iC^i$.
Since $AC^{d-k}$ is similar to a permutation matrix, its trace $d\alpha_k$
is a non-negative integer. Hence $\alpha_k$ for $k=0,1,\cdots,d-1$
are non-negative rationals. By Lemma \ref{permutable},
all the eigenvalues of $A$ are roots of unity, in particular,
$\sum_{i=0}^{d-1}\alpha_i =1 =\big|\sum_{i=0}^{d-1}\alpha_i\lambda^i \big|$.
If there are at least two of the coefficients non-zero, say $\alpha_k\ne 0\ne \alpha_h$
for $0\le k\ne h<d$, then 
$\big|\alpha_k\lambda^k+\alpha_h\lambda^h\big|<\alpha_k+\alpha_h$
because $\lambda^k\ne\lambda^h$; hence
$$
1=\Big|\sum_{i=0}^{d-1}\alpha_i\lambda^i \Big|\le
\big|\alpha_k\lambda^k+\alpha_h\lambda^h\big|+
 \sum_{i\ne k,h}\big|\alpha_i\lambda^i\big|<\sum_{i=1}^{d-1}\alpha_i=1,
$$
which is impossible. Thus there is exactly one of the coefficients, say $\alpha_k$,
which is non-zero; then $\alpha_k=1$ and $A=C^k$.

\section{The case of $p$-groups} \label{case p group}

In the following $p$ is always an odd prime and $n$ is a positive integer.
In this section we consider a special case of Theorem \ref{main},
where the permutation-like matrix group is a $p$-group.
This is a key step for the proof of the theorem.

Let $C\in{\rm GL}_{p^n}({\Bbb C})$ be a maximal cycle.
We have a disjoint union
${\Bbb Z}_{p^n}=p{\Bbb Z}_{p^n}\bigcup {\Bbb Z}_{p^n}^*$, where 
$$p{\Bbb Z}_{p^n}=\{pk~({\rm mod}~p^n)\mid k\in{\Bbb Z}_{p^n}\}
 =\{0,p,\cdots,p(p^{n-1}-1)\}.$$
By Lemma \ref{spectrum}, we have two subspaces of ${\Bbb C}^{p^n}$,
denoted by $V^p$ and $V^*$,  as follows:
\begin{equation}\label{VpV*}
V^{p} =\bigoplus_{j\in p{\Bbb Z}_{p^n}}{\rm E}(\lambda^j)
=\bigoplus_{e\in{\cal E}^p}{\Bbb C}e
\quad\mbox{and}\quad
 V^{*} =\bigoplus_{j\in {\Bbb Z}_{p^n}^*}{\rm E}(\lambda^j)
=\bigoplus_{e\in{\cal E}^*}{\Bbb C}e
\end{equation}
where ${\cal E}^p$ and ${\cal E}^*$ denote the basis
of $V^p$ and $V^*$ respectively as in Lemma \ref{basis}.
Then ${\Bbb C}^{p^n}=V^p\oplus V^*$.

For any $A\in{\rm GL}_{p^n}({\Bbb C})$ which normalizes $\langle C\rangle$,
there is exactly one $r\in{\Bbb Z}_{p^n}^*$ such that $A^{-1}CA=C^r$.
Both $p{\Bbb Z}_{p^n}$ and ${\Bbb Z}_{p^n}^*$ are $\mu_r$-invariant.
So, for any integer $k$ with $0<k<p^n$ both $V^{p}$ and $V^{*}$
are $AC^k$-invariant subspaces of ${\Bbb C}^{p^n}$.
By $AC^k|_{V^p}$ and $AC^k|_{V^*}$ we denote the linear transformations 
of $AC^k$ restricted to $V^p$ and $V^*$ respectively. 
Correspondingly, 
$AC^k|_{{\cal E}^p}$ and $AC^k|_{{\cal E}^*}$ are
matrices of $AC^k|_{V^p}$ and $AC^k|_{V^*}$ respectively.

\begin{lemma} \label{AC^k}
Let $A,C\in{\rm GL}_{p^n}({\Bbb C})$.
Assume that $C$ is a maximal cycle and
$A^{-1}CA=C^r$ for an $r\in{\Bbb Z}_{p^n}^*$
with ${\rm ord}_{{\Bbb Z}_{p^{n}}^*}(r)=p^{a}$ where $0\le a<n$.
Let $V^p$ and $V^*$ be as in Eqn \eqref{VpV*}.
If $A^{p^{a}}=I$,  then for 
$${\rm char}_{(AC^k)|_{V^*}}(x)
=\begin{cases}{\Phi_{p^{n-\nu_p(k)}}(x)}^{p^{\nu_p(k)}}, & 0\le\nu_p(k)<n-a;\\[5pt]
 {(x^{p^a}-1)}^{p^{n-a-1}(p-1)}, &  \nu_p(k)\ge n-a,~{\rm or}~ k=0; \end{cases}
$$
where $\nu_p(k)$ denotes the $p$-adic valuation of $k$, i.e.
$p^{\nu_p(k)}$ is the largest power of~$p$ which divides $k$.
\end{lemma}

\pf Assume that ${\Bbb Z}_{p^n}^*$ is partitioned into $\mu_r$-orbits
$\Gamma_1$, $\cdots$, $\Gamma_h$.
Every orbit $\Gamma_i$ has length $p^{a}$,
and the number $h=p^{n-a-1}(p-1)$. We can assume that
$$
\Gamma_1=\{j_1,rj_1,\cdots,r^{p^{a}-1}j_1\},~  \cdots,~
\Gamma_{h}=\{j_h, rj_h, \cdots, r^{p^{a}-1}j_h\}.
$$
For the basis ${\cal E}^*$ of $V^*$, by Lemma \ref{basis} we have
 ${\cal E}^*={\cal E}^*_1\cup\cdots\cup{\cal E}^*_h$ with
\begin{equation}\label{E_i}
{\cal E}^*_i=\big\{e_{j_i},~Ae_{j_i},~\cdots,~A^{r^{p^a-1}}e_{j_i}\big\},
\qquad i=1,\cdots,h.
\end{equation}
Denote the restricted matrices to ${\cal E}^*_i$ by
$A_i=A|_{{\cal E}^*_i}$ and $C_i=C|_{{\cal E}^*_i}$ for $i=1,\cdots,h$.
Since $A_i^{p^a}=I$, we have
\begin{equation}\label{A_i}
A_{i}=
   \begin{pmatrix}0&\cdots&0&1\\ 1&\ddots&\ddots&0\\
          \vdots&\ddots&\ddots&\vdots\\
          0&\cdots&1&0   \end{pmatrix}_{p^{a}\times p^{a}},~
 C_{i}=
   \begin{pmatrix}\lambda^{j_i} \\ &\lambda^{rj_i}\\
   &&\ddots\\ &&& \lambda^{r^{p^a-1}j_i}
 \end{pmatrix}_{p^{a}\times p^{a}}.
\end{equation}
Thus
$
AC^k|_{{\cal E}^*}=\bigoplus_{i=1}^h A_iC_i^k,
$
and
$$
{\rm char}_{A_iC_i^k}(x)=x^{p^a}-\lambda^{\sum_{j\in\Gamma_i}jk}
=x^{p^a}-\lambda^{j_ik(1+r+\cdots+r^{p^a-1})}.
$$
The conclusion is obviously true if $k=0$. So we further assume that $k\ne 0$.
By Lemma \ref{p-cyclic} (iii), we can take an integer $v$ which is coprime
to~$p$ such that $r=1+vp^{n-a}$. Then it is easy to check that
$\nu_p(1+r+\cdots+r^{p^a-1})=a$, see Eqn \eqref{r-1}.
So we can write $1+r+\cdots+r^{p^a-1}=a'p^a$ and $k=k'p^{\nu_p(k)}$
with $p\nmid a'$ and $p\nmid k'$. Then
$$
{\rm char}_{A_iC_i^k}(x)=x^{p^a}-\lambda^{a'k'j_i p^{\nu_p(k)+a}},\qquad
i=1,\cdots,h.
$$
By Lemma \ref{p-cyclic} (ii),
the collection of $\lambda^{a'k'j_i p^{\nu_p(k)+a}}$ for $i=1,\cdots, h$
is just the collection of all primitive $p^{n-a-\nu_p(k)}$-th roots of unity,
each of which appears with multiplicity
$$\frac{h}{p^{{n-a-\nu_p(k)}-1}(p-1)}=p^{\nu_p(k)}.$$
If $n-a-\nu_p(k)>0$, then the collection of roots of ${\rm char}_{AC^k|_{V^*}}(x)$
is just the collection of all primitive  $p^{n-\nu_p(k)}$-th roots of unity,
each of which has multiplicity $p^{\nu_p(k)}$; hence
$${\rm char}_{AC^k|_{V^*}}(x)={\Phi_{p^{n-\nu_p(k)}}(x)}^{p^{\nu_p(k)}}.$$
Otherwise, $n-a-\nu_p(k)\le 0$, i.e.
$\lambda^{a'k'j_i p^{\nu_p(k)+a}}=1$, hence
$${\rm char}_{AC^k|_{V^*}}(x)=(x^{p^a}-1)^{h}
={(x^{p^a}-1)}^{p^{n-a-1}(p-1)}.\eqno\Box$$

\begin{corollary}\label{V^p}
Let notation be as in Lemma \ref{AC^k}.
If the matrix group $\langle A,C\rangle$ is permutation-like
and $A^{p^a}=I$,  then $\big\langle A|_{V^p},C|_{V^p}\big\rangle$
is a permutation-like matrix group of dimension $p^{n-1}$.
\end{corollary}

\pf Let $\ell,k$ be any non-zero integers.
Let $A'=A^\ell$, $a'=a-\nu_p(\ell)$ and $r'=r^\ell$.
Then $A'^{-1}CA'=C^{r'}$ and
${\rm ord}(A')={\rm ord}_{{\Bbb Z}_{p^n}^*}(r')=p^{a'}$.
by Lemma~\ref{AC^k},
$${\rm char}_{(A^\ell C^k)|_{V^*}}(x)
=\begin{cases}{\Phi_{p^{n-\nu_p(k)}}(x)}^{p^{\nu_p(k)}}, & \nu_p(k)<n-a';\\[5pt]
 {(x^{p^{a'}}-1)}^{p^{n-a'-1}(p-1)}, &  \nu_p(k)\ge n-a'. \end{cases}
$$
Since $A^\ell C^k=A'C^k$ is similar to a permutation matrix,
by Lemma \ref{permutable} and its corollary,
$$
{\rm char}_{A^\ell C^k|_{V^p}}(x)=\begin{cases}
{(x^{p^{n-\nu_p(k)-1}}-1)}^{p^{\nu_p(k)}}, & \nu_p(k)<n-a';\\
\prod_{i}{(x^{p^i}-1)}^{j_i}, & \nu_p(k)\ge n-a'.
\end{cases}
$$
By Lemma \ref{permutable} again, the matrix group
$\big\langle A|_{V^p},C|_{V^p}\big\rangle$
is a permutation-like matrix group of dimension $p^{n-1}$.
\qed

\begin{lemma} \label{A V^p}
Let $A,C\in{\rm GL}_{p^n}({\Bbb C})$.
Assume that $C$ is a maximal cycle and
$A^{-1}CA=C^r$ for an $r\in{\Bbb Z}_{p^n}^*$
with ${\rm ord}_{{\Bbb Z}_{p^{n}}^*}(r)=p$.
If $\langle A,C\rangle$ is a permutation-like matrix group and $A^{p}=I$,
then $A|_{V^p}=I$.
\end{lemma}

\pf Note that $\langle A,C\rangle=\{A^\ell C^k\mid 0\le\ell<p,\;0\le k<p^n\}$.
For any $0<\ell<p$ and $0<k<p^n$, by Lemma \ref{AC^k}
(note that $a=1$ at present case),
$$
{\rm char}_{A^\ell C^k|_{V^*}}(x)=\begin{cases}
{\Phi_{p^{n-\nu_p(k)}}(x)}^{p^{\nu_p(k)}}, & \nu_p(k)<n-1;\\
{(x^p-1)}^{p^{n-2}(p-1)}, & \nu_p(k)\ge n-1.
\end{cases}
$$
Since $A^\ell C^k$ is similar to a permutation matrix, by Lemma \ref{permutable},
\begin{equation}\label{A^p=I}
{\rm char}_{A^\ell C^k|_{V^p}}(x)=\begin{cases}
{(x^{p^{n-\nu_p(k)-1}}-1)}^{p^{\nu_p(k)}}, & \nu_p(k)<n-1;\\
{(x-1)}^{p^{n-1}-pj}{(x^p-1)}^{j}, & \nu_p(k)\ge n-1.
\end{cases}
\end{equation}
Thus $\big\langle A|_{{\cal E}^p},C|_{{\cal E}^p}\big\rangle$
is a permutation-like matrix group of dimension $p^{n-1}$.

Obviously, $C|_{{\cal E}^p}=\bigoplus_{j\in{p{\Bbb Z}_{p^n}}}\lambda^{j}$
is a maximal cycle of dimension $p^{n-1}$.
Since $r=1+vp^{n-1}$ for some integer $v$ coprime to $p$
(see Lemma \ref{p-cyclic}(iii)),
for any $pt\in p{\Bbb Z}_{p^n}$ we have $\mu_r(pt)\equiv pt\pmod{p^n}$.
So, $A|_{{\cal E}^p}$ is a diagonal matrix, hence
$A|_{{\cal E}^p}$ commutes with $C|_{{\cal E}^p}$.
By Lemma \ref{self-centralized},
$A|_{{\cal E}^p}\in\big\langle C|_{{\cal E}^p}\big\rangle$.
But $(A|_{{\cal E}^p})^p=I$ and ${\rm ord}(C|_{{\cal E}^p})=p^{n-1}$.
Thus $A|_{{\cal E}^p}\in\big\langle C_{{\cal E}^p}^{p^{n-2}}\big\rangle$, 
and we can assume that
$A|_{{\cal E}^p}={(C|_{{\cal E}^p})}^{b p^{n-2}}$
with $0\le b<p$. Suppose that $b>0$, then
$\nu_p(-b p^{n-2})=n-2$, and by Eqn \eqref{A^p=I},
$$
{\rm char}_{I_{p^{n-1}\times p^{n-1}}}(x)
={\rm char}_{(A|_{{\cal E}^p})(C|_{{\cal E}^p})^{-b p^{n-2}}}(x)
={(x^p-1)}^{p^{n-2}}.
$$
However, ${\rm char}_{I_{p^{n-1}\times p^{n-1}}}(x)=(x-1)^{p^{n-1}}$.
This is a contradiction. Thus
$b=0$ and $A|_{{\cal E}^p}=I_{p^{n-1}\times p^{n-1}}$.
\qed

\begin{proposition}\label{p-group}
Let $A,C\in{\rm GL}_{p^n}({\Bbb C})$.
Assume that $C$ is a maximal cycle and
$A^{-1}CA=C^r$ for an $r\in{\Bbb Z}_{p^n}^*$
with ${\rm ord}_{{\Bbb Z}_{p^n}^*}(r)=p^a$ where $0\le a<n$.
If ${\cal G}=\langle A, C\rangle$ is a permutation-like matrix group and $A^{p^a}=I$,
then $A|_{\cal E}$ is a permutation matrix.
\end{proposition}

\pf If $n=1$, then $a=0$ and $A=I$, the proposition holds trivially.

Assume $n>1$.
From Eqn \eqref{A_i} we have seen that $A|_{V^*}$ is a permutation matrix.
By Corollary~\ref{V^p}, $\langle A|_{V^p}, C|_{V^p}\rangle$ is a permutation-like
matrix group of dimension $p^{n-1}$. Note that
${\rm ord}_{\cal G}(A^{p^{a-1}})
 ={\rm ord}_{{\Bbb Z}_{p^{n}}^*}(r^{p^{a-1}})=p$.
By Lemma~\ref{A V^p}, $(A|_{V^p})^{p^{a-1}}=I$.
And, ${\rm ord}_{{\Bbb Z}_{p^{n-1}}^*}(r)=p^{a-1}$.
By induction on $n$, $A|_{{\cal E}^p}$ is a permutation matrix.
Hence $A|_{\cal E}=A|_{{\cal E}^p}\oplus A|_{{\cal E}^*}$
is a permutation matrix.
 \qed

\begin{corollary}\label{p group}
Let $A,C\in{\rm GL}_{p^n}({\Bbb C})$.
Assume that $C$ is a maximal cycle and $A$ normalizes $\langle C\rangle$.
If ${\cal G}=\langle A, C\rangle$ is a permutation-like matrix $p$-group,
then ${\cal G}$ is a permutation matrix group.
\end{corollary}

\pf 
By Lemma~\ref{self-centralized}, $\langle C\rangle$ is self-centralized in ${\cal G}$.
By Lemma \ref{split}(i), we have an $A'\in{\cal G}$ such that
${\cal G}=\langle A',C\rangle$, $A'^{-1}CA'=C^r$ for an $r\in{\Bbb Z}_{p^n}^*$
with ${\rm ord}_{{\Bbb Z}_{p^n}^*}(r)=p^a$ and $A'^{p^a}=I$.
Thus, by Proposition \ref{p-group} and Proposition \ref{SC},
${\cal G}$ is a permutation group. \qed

\section{Proof of Theorem \ref{main}} \label{proof}

Let ${\cal G}\le{\rm G}_{p^n}({\Bbb C})$
be a permutation-like matrix group of dimension $p^n$
where $p$ is an odd prime. Let $C\in{\cal G}$ be a maximal cycle
such that $\langle C\rangle$ is a normal  subgroup of ${\cal G}$.
By Lemma \ref{self-centralized}, $\langle C\rangle$ is self-centralized in ${\cal G}$.
By Lemma \ref{split} and Eqn \eqref{r}, ${\cal G}=\langle A,C\rangle$,
$A^{-1}CA=C^r$ where $r\in{\Bbb Z}_{p^n}^*$ with
$$
{\rm ord}_{{\Bbb Z}_{p^n}^*}(r)=|{\cal G}/\langle C\rangle|=sp^a,\qquad
 s\,|\,(p-1),~~0\le a<n.
$$
Let ${\cal E}$ be the basis of ${\Bbb C}^{p^n}$ as in Lemma \ref{basis}.

\medskip{\it Case 1:}~ $s=1$. So ${\cal G}$ is a $p$-group, and
Theorem \ref{main} holds by Corollary \ref{p group}.

\medskip{\it Case 2:}~ $a=0$ and $s>1$.
By Lemma \ref{p-cyclic} (iii), we can assume that
${\rm ord}_{{\Bbb Z}_{p}^*}(r)=s$.
Then for any $0<\ell<s$, $r^\ell\not\equiv 1\pmod{p}$,
i.e. $p\nmid(r^\ell-1)$.
Thus $r^\ell k\equiv k\pmod{p^n}$ if and only if $k\equiv 0\pmod{p^n}$.
So every $\mu_r$-orbit on ${\Bbb Z}_{p^n}$ has length $s$ except for
the orbit $\{0\}$. By Lemma \ref{split} (ii), $A^s=I$.
By Lemma \ref{basis} (ii) and (iii),
$$
A|_{\cal E}=\omega_0\oplus
\begin{pmatrix}0&\cdots&0&1\\ 1&\ddots&\ddots&0\\
    \vdots&\ddots&\ddots&\vdots\\ 0&\cdots&1&0\end{pmatrix}_{s\times s}
\oplus\cdots\oplus
\begin{pmatrix}0&\cdots&0&1\\ 1&\ddots&\ddots&0\\
    \vdots&\ddots&\ddots&\vdots\\ 0&\cdots&1&0\end{pmatrix}_{s\times s}.
$$
And the characteristic polynomial
$${\rm char}_A(x)=(x-\omega_0)(x^s-1)^{(p^n-1)/s}.$$
By Lemma \ref{permutable}, $\omega_0=1$.
So $A|_{\cal E}$ is a permutation matrix.
By Proposition \ref{SC}, we obtain that ${\cal G}$ is a permutation matrix group.

\medskip{\it Case 3:}~ $a>0$ and $s>1$.
Then ${\cal G}$ is not a $p$-group.
By Lemma \ref{split}(ii) we further have that $A^{sp^a}=I$.
Since $s$ and $p^a$ are coprime each other, we have integers
$t,m$ such that $st+p^a m=1$.
Let $A'=A^{st}$ and $A''=A^{p^am}$.
Then $A=A^{st+p^a m}=A'A''$, ${A'}^{p^a}=I$ and ${A''}^s=I$.
By Proposition~\ref{p-group}, $A'|_{\cal E}$ is a permutation matrix.
From the above argument of Case 2, $A''|_{\cal E}$ is also a permutation matrix.
Thus $A|_{\cal E}=A'|_{\cal E}\cdot A''|_{\cal E}$ is a permutation matrix.
By Proposition \ref{SC}, ${\cal G}$ is a permutation matrix group.
\qed

\section*{Acknowledgments}

The authors are supported by NSFC through the grant number 11271005.
Thanks are given to the reviewers for the helpful suggestions.

\end{document}